\documentclass[11pt]{article}
\usepackage{amsmath,amssymb}
\usepackage{amsfonts}
\usepackage{eucal}
\title{\leavevmode\vadjust{\vskip -5mm}
\textbf{Group-theoretic Description of Riemannian Spaces}}
\author{\sc Serhiy E. SAMOKHVALOV\thanks{{\bf e}-{\it mail}: samokhval@dstu.dp.ua}\\
 \tabaddress{Department of Applied Mathematics \\
State Technical University, Dniprodzerzhinsk, Ukraine}}

\date{April 22, 2007}

\def\tabaddress#1{{\small\it\begin{tabular}[t]{c}#1
\\[1.2ex]\end{tabular}}}

\def\<#1>{\langle#1\rangle}

\def\beq{\begin{equation}}
\def\eeq{\end{equation}}
\def\bea{\begin{eqnarray}}
\def\eea{\end{eqnarray}}
\def\beann{\begin{eqnarray*}}
\def\eeann{\end{eqnarray*}}
\def\ben{\begin{enumerate}}
\def\een{\end{enumerate}}

\def\qed{\ifvmode\removelastskip\fi
{\unskip\nobreak\hfil\penalty50\hbox{}\nobreak\hfil \hbox{\vrule
height1.2ex width1.2ex}\parfillskip=0pt \finalhyphendemerits=0
\par\smallskip}}

\def\texthook{\vrule height 0pt depth 0.4pt width 3.5pt
          \vrule height 5pt depth 0.4pt \kern 3pt}
\def\scripthook{\vrule height 0pt depth 0.2pt width 1.5pt
                \vrule height 3pt depth 0.2pt width 0.2pt \kern 1pt}

\begin{document}

\maketitle

\thispagestyle{empty}

\begin{abstract}
\noindent It is shown that a locally geometrical structure of
arbitrarily curved Riemannian space is defined by a deformed group
of its diffeomorphisms.
\end{abstract}


\\

\clearpage


Until recently it was thought impossible to realize Klein's
Erlangen Program \cite{1} for geometrical structures with
arbitrary variable curvature; this Riemann-Klein antagonism, as it
was figuratively called by E. Cartan \cite{2}, could only be
overcome at the cost of program's modification and rejection of
group structure of transformations which were used. Thus in
\cite{3} categories are employed while in \cite{4} quasigroups
are, and it is even stated that quasigroups are an algebraic
equivalent of geometric notion of curvature.

In work~\cite{5} it was shown that group-theoretic description of
connections in fiber bundles with arbitrary variable curvature can
be performed by means of deformed infinite Lie groups introduced
out of physical considerations in work~\cite{6}, the structural
equation follows from group axioms and is a necessary condition
for existence of a group which defines given geometrical
structure. This allowed realization of Klein's Program for
connections in fiber bundles.

The structure of (pseudo)Riemannian space $M$ is a special case of
structure of affine connection in tangent bundle and therefore it
can be specified similarly to arbitrary connection~\cite{5}. At
the same time it necessary to apply additional conditions of
torsion absence and coordination of connection with metric. The
group fulfilling this description acts in tangent bundle of space
$M$ and is an infinite and specially deformed group which has the
structure of semidirect product of diffeomorphisms group $\Gamma_T
= Diff\ M$ and gauge group $SO \left(m, n-m \right)^g$, where $n$
is space dimension~\cite{6}.

It was shown in~\cite{6} that there exists a more natural way of
group-theoretic description for (pseudo)Riemannian spaces, the one
with the help of a narrower group, i.e. the deformed group
$\Gamma^H_T$ of diffeomorphisms of space $M$. The generators of
such group define on M (locally, within the bounds of coordinate
chart) field of affine vielbein, multiplication law define the
rule of parallel transport of vectors, in consideration of which
torsion is automatically zeroed in view of group axioms,
components of vielbeins field in coordinate basis as well as
anholonomity and connection coefficients are expressed through
auxiliary deformation functions by means of which group
$\Gamma^H_T$ is built. Locally any space of torsion-free affine
connection can be described in such fashion. With additional
assumption that vielbeins field, defined by action of group
$\Gamma^H_T$ is (pseudo)orthonormal, and in case of parallel
transport of vectors they merely rotate, coefficients of affine
connection in coordinate basis automatically become Christoffel
symbols, i.e. they are defined through metric in a certain way,
therefore there is no need to postulate this statement.

Publication~\cite{6} had physical value and group-theoretical and
geometrical aspects were only slightly touched upon there, while
some important geometrical relations were neglected at all. This
work makes up for this. Specifically, we show that definition of
curvature tensor (which in our approach becomes a characteristic
of group $\Gamma^H_T$) through connection coefficients follows
from equation which comes from group axioms and is an essential
condition for existence of group $\Gamma^H_T$.

With the help of groups $\Gamma^H_T$ Klein's Erlangen Program is
realized for (pseudo)Riemannian spaces of arbitrary variable
curvature, in the most rational fashion at that. Groups
$\Gamma^H_T$ act on $M$ and their transformations are interpreted
as gauge translation in curved (pseudo)Riemannian spaces. It is
due to this that group-theoretic description of (pseudo)Riemannian
spaces through groups $\Gamma^H_T$ is important for gravitation
theory, gravitation being interpreted as gauge theory of
translations group~\cite{7}.

The work doesn't deal with global topological problems and all
relations are obtained within the bounds of a single coordinate
chart. Besides, we perceive groups to be respective local groups.

\textbf{1.} Let's specify the general procedure of building
deformed infinite Lie groups~\cite{5} for the case of deformed
group of diffeomorphisms $\Gamma^H_T$. This time, opposite
to~\cite{5} we will use coordinate approach.

Let $O$ be a coordinate chart on manifold $M$ with coordinates
$x^{\mu}$ (we use Greek alphabet for indices). We assume
coordinates to be fixed and won't change them further.

In $O$ there acts Abelian group of translations $T=\{\tilde{t} \}$
 according to the formula:
\[x'^{\mu}=x^{\mu}+\tilde{t}^{\mu}
\]
\noindent In set $C_{\infty}(O, T)$ of smooth mappings of $O$ in
$T$ let's single out subset $\Gamma_T = \{\tilde{t}(x)\}$ with
condition:
\[det\{{\delta}^{\mu}_{\nu} + \partial_{\nu}
\tilde{t}^{\mu} (x)\} \ne 0,\ \forall x \in O,
\]
\noindent where $\partial_\nu :=\partial /
\partial x^\nu$, and assign to it
the multiplication law $\tilde{t''} = \tilde{t} \times
\tilde{t'}$:
\begin{equation}\label{eq1}
\tilde{t''}^\mu (x) = \tilde{t}^\mu (x) + \tilde{t'}^\mu (x'),
\end{equation}
\noindent where
\begin{equation}\label{eq2}
{x'}^\mu = {x}^\mu + \tilde{t}^\mu (x).
\end{equation}
\noindent With it the set $\Gamma_T$ becomes a local group. Group
$\Gamma_T$ acts smoothly in chart $O$ according to the formula
(\ref{eq2}) and is a local group of diffeomorphisms of chart $O$
in additive parameterization. According to definition 1 from [5]
group $\Gamma_T$ is \emph{a group of undeformed chart $O$, or
undeformed group}.

Let's deform group $\Gamma_T$ by means of \emph{deformation} $H$
which is defined by mapping $H : O \times T \rightarrow T$ with
properties which in our case are described as follows:

$1H)\ H \in C_{\infty}(O, T)$;

$2H)\ H(x,0)=0,\ \forall x \in O$;

$3H)\ \exists\ {\rm mapping}\ K : O \times T \to T: K(x, H(x,
\tilde{t})) = \tilde{t},\ \forall x \in O,\ \tilde{t} \in T.$

Group $\Gamma^H_T=\{t(x)\}$ is obtained from group $\Gamma_T$ by
isomorphism, which is specified by deformation $H$ according to
the formula:
\begin{equation}\label{eq3}
t^m(x)=H^m(x,\tilde t (x)).
\end{equation}
\noindent Functions $t^m(x)$ which parameterize group $\Gamma^H_T$
(we use indices from Latin alphabet for them), satisfy the
condition:
\[det\{\delta^{\mu}_{\nu}+d_{\nu}K^{\mu}(x, t(x))\}\ne 0,\ \forall
x \in O,
\]
\noindent where $d_{\nu} :=d/dx^{\nu}$, and multiplication law
$t'' = t*t'$ is defined by isomorphism (\ref{eq3}):
\begin{equation}\label{eq4}
t''^m(x)= \varphi^m (x,t(x),t'(x')):= H^m(x,K(x,t(x))+K(x',t'
(x'))),
\end{equation}
\noindent where
\begin{equation}\label{eq5}
{x'}^{\mu} = f^{\mu}(x,t(x)):= x^{\mu} + K^{\mu}(x,t(x)).
\end{equation}
\noindent Group $\Gamma^H_T$ acts smoothly in the chart $O$
according to the formula (\ref{eq5}). According to definition 3
from~\cite{5} group $\Gamma^H_T$ is \emph {a group of deformed
chart} $O$, or \emph {deformed group}.

Multiplication law (\ref{eq4}) for deformed group $\Gamma^H_T$
explicitly depends on $x$, and, therefore, structural constants
analogue for groups $\Gamma^H_T$ is structure functions
$F(x)^n_{kl}$, which are defined by the formula:
\begin{equation}\label{eq6}
F(x)^n_{kl}:=
\left(\partial^2_{k,l'}-\partial^2_{l,k'}\right)\varphi^n(x,t,t')\Bigr|_{t=t'=0}.
\end{equation}
\noindent (here and henceforth $\partial_k:=\partial / \partial
t^k$, primed index stands for differentiation with respect to
$t'$).

\textbf{2.} Let's introduce auxiliary functions:
\[h(x)^m_{\mu} =
\frac{\displaystyle \partial}{\displaystyle
\partial \tilde t^{\mu}}H^m(x,\tilde t)\Bigr |_{\tilde t=0}.
\]
\noindent Property $3H$ allows fulfillment of condition:
\begin{equation}\label{eq7}
det \left \{h(x)^m_{\mu}\right \}\ne 0,\ \ \forall x \in O,
\end{equation}
\noindent wherefrom there follows existence of functions
$h(x)^{\mu}_m$ of the type that
$h(x)^{\mu}_nh(x)_{\mu}^m=\delta^m_n$, $\forall x \in O$. It is
obvious that $h(x)^{\mu}_m=\partial_m K^{\mu}(x, t)\Bigr |_{t=0}$.
With the help of these functions we will substitute Greek indices
for Latin and vice versa.

Assuming parameters in multiplication law for group $\Gamma^H_T$
to be constant, let's define functions:
\begin{equation}\label{eq8}
\mu(x, t)^m{}_n := \partial_{n'} \varphi^m(x,t',t) \Bigr |_{t'=0},
\end{equation}
\begin{equation}\label{eq9}
\lambda(x, t)^m{}_n := \partial_{n'} \varphi^m(x,t,t') \Bigr
|_{t'=0}.
\end{equation}

The condition of multiplication law associativity in group
$\Gamma^H_T$: $(t*t')*t''=t*(t'*t'')$ is fulfilled automatically
for any deformation H in view of multiplication law (\ref{eq1})
associativity in diffeomorphisms group. Let's derive this
condition for constant parameters of group $\Gamma^H_T$:
\begin{equation}\label{eq10}
\varphi^m(x,\varphi (x,t,t'),t'')=\varphi^m(x,t,\varphi
(x',t',t'')).
\end{equation}
\noindent Differentiating it with respect to t in zero we obtain
the equation:
\begin{equation}\label{eq11}
h(x)^{\mu}_k \partial_{\mu} \varphi^m(x,t,t')-\mu(x,t)^nk
\partial_n\varphi^m(x,t,t')=-\mu (x,\varphi(x,t,t'))^m{}_k,
\end{equation}
\noindent and with respect to $t''$ in zero the equation:
\begin{equation}\label{eq12}
\lambda(x', t')^n{}_k \partial_{n'} \varphi^m(x,t,t')=\lambda
(x,\varphi (x,t,t'))^m{}_k.
\end{equation}

The condition for their integrabilily is equation:
\[h(x)^{\nu}_k \partial_{\nu}\mu (x,t)^m{}_l-\mu(x,t)^n{}_k\partial_n
\mu(x,t)^m{}_l - h(x)^{\nu}_l \partial_{\nu} \mu(x,t)^m{}_k+\mu
(x,t)^n{}_l \partial_n \mu (x,t)^m{}_k=
\]
\begin{equation}\label{eq13}
=F(x)^n_{kl}\mu(x,t)^m{}_n
\end{equation}
\noindent and
\begin{equation}\label{eq14}
\lambda(x, t)^n{}_k \partial_{n}\lambda(x,t)^m{}_l-
\lambda(x,t)^n{}_l\partial_n
\lambda(x,t)^m{}_k=F(x')^n_{kl}\lambda(x,t)^m{}_n.
\end{equation}
\noindent respectively.

Let's call equations (\ref{eq11}) and (\ref{eq12}) \emph{the left
and the right Lie equation for groups} $\Gamma^H_T$, while
equations (\ref{eq13}) and (\ref{eq14}) \emph{the left and the
right Maurer-Cartan equations for groups} $\Gamma^H_T$.

If the condition of associativity (\ref{eq10}) is immediately
differentiated with respect to $t'$ and $t''$ in zero with
differing sequence we obtain the equation:
\begin{equation}\label{eq15}
h(x)^{\nu}_k \partial_{\nu}\lambda
(x,t)^m{}_l-\mu(x,t)^n{}_k\partial_n \lambda(x,t)^m{}_l +
\lambda(x,t)^{n}{}_l
\partial_{n} \mu(x,t)^m{}_k=0.
\end{equation}

Let's perform consequently two $\Gamma^H_T$-transformations with
constant parameters $t$ and $t'$. Composition law of
transformations results in equation:
\[f^\mu (f(x,t),t')=f^\mu(x,\varphi(x,t,t')),
\]
\noindent which is fulfilled automatically for any deformation $H$
in view of performance of composition law in the group of
diffeomorphisms $\Gamma_T$. Differentiating it with respect to $t$
in zero we obtain the equation:
\begin{equation}\label{eq16}
h(x)^{\nu}_k \partial_{\nu}f^\mu (x,t)-\mu(x,t)^n{}_k\partial_n
f^\mu(x,t)=0,
\end{equation}
\noindent and differentiating it with respect to $t'$ in zero the
equation:
\begin{equation}\label{eq17}
h(x')^{\mu}_k -\lambda (x,t)^n{}_k\partial_n f^\mu(x,t)=0.
\end{equation}
\noindent The condition of integrabilily of these equations in
case of fulfillment of equations (\ref{eq13}) and (\ref{eq14}) is
equation:
\begin{equation}\label{eq18}
h(x)^{\nu}_k \partial_\nu h(x)^\mu_l-h(x)^{\nu}_l \partial_\nu
h(x)^\mu_k =F(x)^{n}_{kl}h(x)^\mu_n.
\end{equation}

We will call equations (\ref{eq16}) and (\ref{eq17}) \emph{the
left and the right Lie equations for groups $\Gamma^H_T$
transformations}, while equation (\ref{eq18}) \emph{the
Maurer-Cartan equation for groups $\Gamma^H_T$ transformations}.

\textbf{3.} Let's introduce differentiating operators:
\[X^\tau_k=h(x)^{\nu}_k \partial_{\nu}-\mu(x,t)^n{}_k\partial_n,
\]
\[X^\upsilon_k=\lambda(x,t)^{n}{}_k \partial_{n},
\]
\noindent which we will call \emph{generators of leftward and
rightward shifts, or horizontal and vertical generators} of group
$\Gamma^H_T$ respectively, as well as
\[X_k=h(x)^{\nu}_k \partial_{\nu}
\]
\noindent - \emph{generators of action of group $\Gamma^H_T$ on}
$O$. In terms of generators, equations (\ref{eq13}) - (\ref{eq15})
as well as equation (\ref{eq18}) have quite an elegant form:
\begin{equation}\label{eq19}
\left[X^\tau_k,X^\tau_l\right]=F(x)^n_{kl}X^\tau_n,
\end{equation}
\begin{equation}\label{eq20}
\left[X^\upsilon_k,X^\upsilon_l\right]=F(x')^n_{kl}X^\upsilon_n,
\end{equation}
\begin{equation}\label{eq21}
\left[X^\tau_k,X^\upsilon_l\right]=0,
\end{equation}
\begin{equation}\label{eq22}
\left[X_k,X_l\right]=F(x)^n_{kl}X_n,
\end{equation}
\noindent where square brackets stand for operators commutator.
These equations follow from multiplication law associativity for
group $\Gamma^H_T$, however, due to its infinity, generators
commutators are expanded into generators not by means of structure
constants as in finite parametric Lie groups, but by means of
structure functions dependent on $x$.

The condition for integrability of equations (\ref{eq19}) -
(\ref{eq22}) is the equation for structure functions of group
$\Gamma^H_T$:
\begin{equation}\label{eq23}
h(x)^{\nu}_k
\partial_\nu F(x)^n_{lm}+F(x)^n_{kp}F(x)^p_{lm}+\mbox{cycle}(klm)=0,
\end{equation}
\noindent which is derived from Jacobi's identity for dual
generators commutator.

\textbf{4.} Let's study the expansion of functions defined by
formulae (\ref{eq8}), (\ref{eq9}) according to group parameters
with accuracy to the second order inclusive:
\begin{equation}\label{eq24}
\mu(x,t)^m{}_n=\delta^m_n+\gamma^m{}_{nk}t^k+\frac{1}{2}
\rho^m{}_{lkn}t^lt^k,
\end{equation}
\begin{equation}\label{eq25}
\lambda(x,t)^m{}_n=\delta^m_n+\gamma^m{}_{kn}t^k+\frac{1}{2}
\sigma^m{}_{lkn}t^lt^k.
\end{equation}
\noindent The coefficients of these expansions $\gamma^m{}_{nk}$,
$\rho^m{}_{nkl}$ and $\sigma^m{}_{nkl}$ depend on $x$ in general
case; however, this dependence, defined by deformation functions,
will be specified a little later and to make it shorter we will
not show explicitly neither to the coefficients themselves nor to
the functions they define. On inserting expansions (\ref{eq24})
and (\ref{eq25}) into formulae (\ref{eq13}) and (\ref{eq14}) we
arrive at the result that in zeroth order with respect to $t$ the
structure functions of group $\Gamma^H_T$ are defined by
skew-symmetric part of coefficients $\gamma^m{}_{kn}$:
\begin{equation}\label{eq26}
F^m_{kn}=\gamma^m{}_{kn}-\gamma^m{}_{nk}.
\end{equation}
\noindent This formula follows directly from definition
(\ref{eq6}) for structure functions of group  $\Gamma^H_T$ and
actually can be considered their definition. Let's introduce
functions:
\begin{equation}\label{eq27}
R^m{}_{lkn}:=\rho^m{}_{lkn}-\rho^m{}_{lnk},
\end{equation}
\begin{equation}\label{eq28}
S^m{}_{lkn}:=\sigma^m{}_{lkn}-\sigma^m{}_{lnk},
\end{equation}
\noindent which we will call \emph{tensors of left and right
curvature of group} $\Gamma^H_T$ respectively. In the first order
with respect to $t$ from formula (\ref{eq13}) we derive:
\begin{equation}\label{eq29}
R^m{}_{lkn}=-\gamma^m{}_{sl}F^s_{kn}+h^\sigma_k\partial_\sigma
\gamma^m{}_{nl}-h^\sigma_n\partial_\sigma\gamma^m{}_{kl}+
\gamma^m{}_{ks}\gamma^s{}_{nl}-\gamma^m{}_{ns}\gamma^s{}_{kl},
\end{equation}
\noindent and from formula (\ref{eq14})
\begin{equation}\label{eq30}
S^m{}_{lkn}=\gamma^m{}_{ls}F^s_{kn}+h^\sigma_l\partial_\sigma
F^m_{kn}+\gamma^m{}_{sk}\gamma^s{}_{ln}-\gamma^m{}_{sn}\gamma^s{}_{lk}.
\end{equation}
\noindent Relation (\ref{eq15}) yields:
\[\sigma^m{}_{lkn}-\rho^m{}_{lnk}=h^\sigma_k\partial_\sigma
\gamma^m{}_{ln}+\gamma^m{}_{ks}\gamma^s{}_{ln}-\gamma^m{}_{sn}\gamma^s{}_{kl},
\]
\noindent wherefrom follows:
\begin{equation}\label{eq31}
R^m{}_{lkn}+S^m{}_{lkn}=h^\sigma_k\partial_\sigma
\gamma^m{}_{ln}-h^\sigma_n\partial_\sigma
\gamma^m{}_{lk}+\gamma^m{}_{ks}\gamma^s{}_{ln}-\gamma^m{}_{sn}
\gamma^s{}_{kl}+\gamma^m{}_{sk}\gamma^s{}_{nl}-\gamma^m{}_{ns}\gamma^s{}_{lk}.
\end{equation}
\noindent Taking into account formulae (\ref{eq26}), (\ref{eq29})
and (\ref{eq30}), expression (\ref{eq31}) is the result of
condition (\ref{eq23}).

\textbf{5.} The relations derived so far, follow solely from group
axioms without consideration of deformation mode of building the
group $\Gamma^H_T$ which allows their fulfillment. However, both
multiplication law (\ref{eq4}) and action (\ref{eq5}) of deformed
group $\Gamma^H_T$ in chart $O$ are defined by deformation $H$
with the help of which it is built. Let's express auxiliary
functions of group $\Gamma^H_T$ through deformation functions. To
this end, let's introduce matrices $H(x,t)^m_{\mu} =
\frac{\displaystyle \partial}{\displaystyle \partial \tilde
t^{\mu}}H^m(x,\tilde t)\Bigr |_{\tilde t=K(x,t)}$. Matrices
$H(x,t)^{\mu}_m =\partial_m K^{\mu}(x,t)$ will be inverse to them.
Direct use of the second equality in (4) in definitions (8) and
(9) yields:
\begin{equation}\label{eq32}
\mu(x,t)^m{}_n=H(x,t)^m_\mu(\delta^\mu_\nu+\partial_\nu
K^\mu(x,t))h(x)^\nu_n,
\end{equation}
\begin{equation}\label{eq33}
\lambda(x,t)^m{}_n=H(x,t)^m_\mu h(x+K(x,t))^\mu_n,
\end{equation}
\noindent or depending upon $\tilde t$:
\[\mu(x,\tilde t)^m{}_n=\frac{\partial}{\partial\tilde
t^\mu}H^m(x,\tilde t)(\delta^\mu_\nu + \partial_\nu \tilde
t^\mu)h(x)^\nu_n,
\]
\begin{equation}\label{eq34}
\lambda(x,\tilde t)^m{}_n=\frac{\partial}{\partial\tilde
t^\mu}H^m(x,\tilde t)h(x+\tilde t)^\mu_n.
\end{equation}

Let's consider the expansion of functions of deformation $H$ up to
the third order with respect to $\tilde t$ inclusive:
\begin{equation}\label{eq35}
H^m(x,\tilde t)=h^m_\mu(\tilde t^\mu+\frac{1}{2}\Gamma^\mu_{\nu
\rho}\tilde t^\nu\tilde t^\rho +\frac{1}{6}\Delta^\mu_{\nu \rho
\sigma}\tilde t^\nu \tilde t^\rho \tilde t^\sigma).
\end{equation}
\noindent Coefficients $h^m_\mu$ satisfy condition (\ref{eq7}) and
$\Gamma^\mu_{\nu\rho}$, $\Delta^\nu_{\nu\rho\sigma}$ are symmetric
in lower indices. On fulfilling these conditions, the coefficients
of expansion (\ref{eq35}) are arbitrary smooth functions of $x$.
Applying them, with accuracy to the second order with respect to
$t$ we derive:
\[K^\mu(x,t)=h^\mu_k t^k-\frac{1}{2}\Gamma^\mu_{kl}t^k t^l,
\]
\[H(x,t)^m_\mu=h^m_\mu +\Gamma^m_{\mu k}t^k
+\frac{1}{2}(\Delta^m_{\mu kl}-\Gamma^m_{\mu s}\Gamma^s_{kl})t^k
t^l.
\]
\noindent In consideration of these expansions formulae
(\ref{eq32}) and (\ref{eq33}) give the following expressions for
coefficients of expansions (\ref{eq24}) and (\ref{eq25}) through
coefficients of expansion (\ref{eq35}):
\begin{equation}\label{eq36}
\gamma^m{}_{kn}=h^m_\mu (\Gamma^\mu_{kn}+h^\nu_k \partial_\nu
h^\mu_n),
\end{equation}
\begin{equation}\label{eq37}
\rho^m{}_{lkn}=h^m_\mu
(\Delta^\mu_{lkn}-\Gamma^\mu_{ns}\Gamma^s_{kl}-h^\nu_n
\partial_\nu \Gamma^\mu_{\kappa\lambda}h^\kappa_k h^\lambda_l),
\end{equation}
\[\sigma^m{}_{lkn}=h^m_\mu
(\Delta^\mu_{lkn}-\Gamma^\mu_{ns}\Gamma^s_{kl}+h^\kappa_k
h^\lambda_l\partial^2_{\kappa\lambda}h^\mu_n-\Gamma^\nu_{kl}
\partial_\nu h^\mu_n + (\Gamma^\mu_{k\sigma}h^\nu_l +\Gamma^\mu_{l\sigma}
h^\nu_k)\partial_\nu h^\sigma_n).
\]
\noindent Inserting these expressions into definitions
(\ref{eq26}) - (\ref{eq28}) and taking into account the symmetry
in lower indices of coefficients $\Gamma^\mu_{kn}$ and
$\Delta^\mu_{lkn}$ we derive formula (18) for structure functions
of group $\Gamma^H_T$, and for its curvature tensors the formulae
as follows:
\begin{equation}\label{eq38}
R^m{}_{lkn}=h^m_\mu
(\partial_\kappa\Gamma^\mu_{\nu\lambda}-\partial_\nu\Gamma^\mu_{\kappa\lambda}+
\Gamma^\mu_{\kappa\sigma}\Gamma^\sigma_{\nu\lambda}-
\Gamma^\mu_{\nu\sigma}\Gamma^\sigma_{\kappa\lambda})h^\lambda_l
h^\kappa_k h^\nu_n,
\end{equation}
\[
S^m{}_{lkn}=h^m_\mu
(\Gamma^\mu_{\l\sigma}F^{\sigma}_{kn}+h^{\sigma}_l
\partial_{\sigma}F^\mu_{kn}+
\Gamma^\mu_{k \sigma}\Gamma^\sigma_{nl}-\Gamma^\mu_{n
\sigma}\Gamma^\sigma_{kl}+ \Gamma^\sigma_{nl}\partial_\sigma
h^{\mu}_k- \Gamma^\sigma_{kl}\partial_\sigma h^{\mu}_n+
\]
\[
h^{\sigma}_l (\Gamma^\mu_{k \sigma} \partial_{\lambda}
h^{\sigma}_n- \Gamma^\mu_{n \sigma}\partial_\lambda
h^{\sigma}_k+\partial_\lambda h^{\sigma}_n \partial_\sigma
h^{\mu}_k-\partial_\lambda h^{\sigma}_k \partial_\sigma
h^{\mu}_n)),
\]
\noindent These formulae could be derived directly from formulae
(\ref{eq29}), (\ref{eq30}) on inserting expressions (\ref{eq36})
and (18) into them. The reason for this is that the condition of
multiplication law associativity in groups $\Gamma^H_T$, which
yields equations (\ref{eq13}) and (\ref{eq14}) wherefrom formulae
(\ref{eq29}) and (30) were derived, is fulfilled automatically
through deformation mode of building groups $\Gamma^H_T$ which we
apply.

\textbf{6.} The very form of tensor of left curvature of group
$\Gamma^H_T$ (formulae (\ref{eq29}) or (\ref{eq38})) as well as
that of its other characteristics indicates that groups
$\Gamma^H_T$ possess ample geometric data which we proceed to
study below.

Generators $X_m=h^\mu_m\partial_\mu$ of action of group
$\Gamma^H_T$ specify on $O$ a field of affine vielbeins, auxiliary
functions of deformation $h^\mu_m$ transfer from coordinate to
affine bases. Elements $t$ of group $\Gamma^H_T$ specify on $O$
\emph{vector fields}  $t=t^mX_m$, parameters $t^m$ of group
$\Gamma^H_T$ are components of these fields in basis $X_m$.

Structure functions of group $\Gamma^H_T$  with lower coordinate
indices in view of formula (\ref{eq18}) can be represented as:
\[F^k_{\mu\nu}=\partial_\nu h^k_\mu -\partial_\mu h^k_\nu,
\]
thus they have geometric meaning (with accuracy to factor -2) of
\emph{anholonomity object}.

Let's study multiplication law  $t*\tau$ in group $\Gamma^H_T$ for
the case of infinitesimal second multiplier:
\[(t*\tau)^m(x)=t^m(x)+\lambda(x,t(x))^m{}_n \tau^n (x'),
\]
\noindent where $x'^\mu=f^\mu(x,t(x))$. Thus, this law gives the
rule for composition of vectors fitted in different points, or the
\emph{rule of parallel transport} of vector field $\tau$ from
point $x'$ to point $x$:
\begin{equation}\label{eq39}
\tau^m_\| (x)=\lambda(x,t(x))^m{}_n \tau^n(x').
\end{equation}
\noindent Taking $t$ to be infinitesimal as well, and taking into
account expansion (\ref{eq25}) we have:
\begin{equation}\label{eq40}
\tau^m_\| (x)=\tau^m(x)+t^n (x) \nabla_n \tau^m (x),
\end{equation}
\noindent where
\[\nabla_n \tau^m (x)=h^\sigma_n
\partial_\sigma \tau^m (x) +\gamma^m{}_{nk} \tau^k (x)
\]
\noindent by definition is a \emph{covariant derivative} of vector
field $\tau$ to the direction $X_n$. Thus, functions
$\gamma^m{}_{nk}$ which define the second of parameters order of
multiplication law in group $\Gamma^H_T$ get geometric meaning of
\emph{coefficients of affine connection in basis} $X_n$.

In coordinate basis, relation (\ref{eq39}) in view of (\ref{eq34})
becomes:
\[\tau^\mu_\| (x)=\frac{\partial}{\partial\tilde t^\nu}H^\mu
(x,\tilde t)\tau^\nu (x+\tilde t),
\]
\noindent or in case of infinitesimal $\tilde t$:
\[\tau^\mu_\| (x)=\tau^\mu (x) +\tilde t^\nu (x) \nabla_\nu \tau^\mu (x),
\]
\noindent in relation to which
\[\nabla_\nu \tau^\mu (x)=\partial_\nu \tau^\mu (x)+\Gamma^\mu_{\sigma\nu}
\tau^\sigma (x).
\]
\noindent Thus, coefficients $\Gamma^\mu_{\sigma\nu}$ which define
the second order of expansion (\ref{eq35}) of deformation
functions get geometric meaning of \emph{coefficients of affine
connection in coordinate basis}. They are arbitrary smooth
functions symmetric in lower indices, corresponding to arbitrary
torsion-free affine connection. For undeformed group $\Gamma^H_T$
covariant derivatives is obviously congruent with partial
derivatives.

It is in this meaning that, specifying the rule of parallel
transport of vectors by its multiplication law (which is defined
by deformation $H$), deformed groups of diffeomorphisms
$\Gamma^H_T$ specify by their action a structure of torsion-free
affine connection in tangent bundle of chart $O$; arbitrary
torsion-free affine connection can be specified over $O$ in such
fashion.

The same connection is specified by all groups $\Gamma^{H'}_T$ in
which coefficients $\Gamma^\mu_{\sigma\nu}$ in the second order of
expansion (\ref{eq35}) of their function of deformation $H'$ are
congruent, particularly if
\begin{equation}\label{eq41}
H'^m(x,\tilde t)=L(x)^m{}_n H^n (x,\tilde t).
\end{equation}
\noindent where matrices $L(x)^m{}_n$, dependent upon $x$, belong
to gauge group $GL(n)^g$. In transformation (\ref{eq41}) the
affine vielbeins field on $O$ changes: $X'_m =L^{-1} (x)^n{}_m
X_n$. The third and higher orders of parameter expansion of
deformation functions do not influence the connection and can be
arbitrary. This is related to the fact that definition
(\ref{eq39}) allows to make parallel transport of vector field
from point $x'$ to point $x$ for \emph{finite} distance $\tilde t
(x) = x'-x=K(x,t(x))$, though infinitesimal shifts are enough to
specify a connection. There are, however, quite natural additional
requirements to deformation functions which follow from geometric
point of view and allow to completely fix deformation functions
with respect to the first two orders of expansion (\ref{eq35}),
i.e. with respect to the affine vielbeins field and affine
connection coefficients. They are related to the generation of
finite parallel transports (\ref{eq39}) with the help of integral
sequence of infinitesimal transports (\ref{eq40}), and we make
plans to study this problem for the structure of affine connection
in our next work \cite{8}.

Let's choose points  $x_1=x+\tilde t_1$,  $x_2=x+\tilde t_2$,
$x_3=x+\tilde t_1+\tilde t_2$ ($\tilde t_1$ and $\tilde t_2$  we
assume to be constant) and perform, according to formula
(\ref{eq39}), parallel transport of vector field $\tau (x)$ from
point $x_3$ to point $x_1$, and then to point $x$ (first choice),
as well as from point $x_3$ to point $x_2$ and then to point $x$
(second choice). The difference of the results obtained gives:
\[\tau^m_\| (x)_1-\tau^m_\|
(x)_2=(\lambda(x,\tilde t_1)^m{}_k \lambda(x_1,\tilde
t_2)^k{}_n-\lambda(x,\tilde t_2)^m{}_k \lambda(x_2,\tilde
t_1)^k{}_n)\tau^n(x_3).
\]
\noindent For infinitesimal  $\tilde t_1$ and $\tilde t_2$ using
formulae (\ref{eq34}), (\ref{eq35}) and (\ref{eq38}) we derive:
\[\tau^m_\| (x)_1-\tau^m_\|
(x)_2=R^m{}_{n\rho\sigma}\tau^n(x)\tilde t_1^\rho \tilde
t_2^\sigma.
\]
\noindent Thus, the tensor  $R^m{}_{n\rho\sigma}$ of the left
curvature of group  $\Gamma^H_T$, which according to the formula
(\ref{eq27}) is a skew-symmetric part of coefficients
$\rho^m{}_{n\rho\sigma}$, which (partially) define the third of
parameter order of multiplication law in group  $\Gamma^H_T$,
acquires the geometric meaning of \emph{curvature tensor} of
affine connection structure, which is specified in $O$ by the
action of group  $\Gamma^H_T$.

Let's summarize the obtained results.

\textbf{Theorem 1.} \emph{Deformed group $\Gamma^H_T$ of
diffeomorphisms in chart $O$ specifies by its action on $O$ an
affine vielbeins field and structure of torsion-free affine
connection in tangent bundle over $O$. Geometric characteristics
of space $O$, such as anholonomity object, affine connection
coefficients, curvature tensor are defined by multiplication law
in group $\Gamma^H_T$, which, in its turn, is defined by
deformation $H$, with the help of which group $\Gamma^H_T$  is
built.}

\emph{Arbitrary torsion-free affine connection can be specified
over $O$ in such fashion.}

Thus, geometric structure of torsion-free affine connection with
arbitrary variable curvature can be referred to only in terms of
deformed groups $\Gamma^H_T$  of diffeomorphisms, due to which
Klein's Erlangen Program is realized for such structure, the
condition of torsion absence (\ref{eq26}) is fulfilled in view of
group axioms and there is no need for its additional application.

\textbf{7.} Let's now assume that matrices $\lambda(x,t)^m{}_n$
belong to the gauge group $SO(m,n-m)^g$, so they satisfy the
equation:
\begin{equation}\label{eq42}
\lambda(x,t)^k{}_m\lambda(x,t)^l{}_n \eta_{kl}=\eta_{mn},
\end{equation}
\noindent where $\eta_{mn}$  is a flat metric (with the help of
which we will lowering indices). This means that vielbein field
$X_m$, specified by the action of group  $\Gamma^H_T$ is
(pseudo)orthonormal and in case of parallel transport of vectors
(\ref{eq39}) they merely (pseudo)rotate. Thus, the action of group
$\Gamma^H_T$ in $O$ specifies the structure of (pseudo)Riemann
space with metric $g_{\mu\nu}=h^m_\mu h^n_\nu \eta_{mn}$.

In the first order with respect to $t$ the equation (\ref{eq42})
produces:
\[\gamma^\centerdot_{ksl}+\gamma^\centerdot_{lsk}=0,
\]
\noindent which allows, with the use of definition (\ref{eq26}),
to express affine connection coefficients in vielbein basis in
terms of structure functions of group  $\Gamma^H_T$:
\begin{equation}\label{eq43}
\gamma^\centerdot_{slk}=\frac{1}{2}\left(F^\centerdot_{slk}+
F^\centerdot_{ksl}+F^\centerdot_{lsk}\right).
\end{equation}
\noindent Recalling geometric interpretation of structure
functions we can see that coefficients $\gamma^s{}_{lk}$  in this
case become \emph{Ricci rotation coefficients}.

With the use of formula (\ref{eq34}) equation (\ref{eq42}) becomes
equation directly for deformation functions:
\begin{equation}\label{eq44}
\frac{\partial}{\partial\tilde t^\mu} H^m (x,\tilde t)
\frac{\partial}{\partial\tilde t^\nu}H^n(x,\tilde
t)\eta_{mn}=g(x+\tilde t)_{\mu\nu}.
\end{equation}
\noindent Besides, we have relation $2H$ for the function $H$:
\begin{equation}\label{eq45}
H^m(x,0)=0,
\end{equation}
\noindent which we will consider a boundary condition for
differential equation (\ref{eq44}). The solution to the problem
(\ref{eq44}), (\ref{eq45}) allows to find deformation functions
$H^m(x,\tilde t)$ with whose help group $\Gamma^H_T$ is produced;
the group specifies in $O$ a structure of (pseudo)Riemannian
space, arbitrary (pseudo)Riemannian structure can be specified in
$O$ in such fashion.

Equation (\ref{eq44}) is invariant under given transformations:
\begin{equation}\label{eq46}
H'^m(x,\tilde t)=\Lambda(x)^m{}_nH^n(x,\tilde t),
\end{equation}
\noindent where matrices  $\Lambda(x)^m{}_n$, dependent upon $x$,
belong to gauge group $SO(m,n-m)^g$, i.e. satisfy relation
$\Lambda(x)^k{}_m\Lambda(x)^l{}_n\eta_{kl}=\eta_{mn}$. Thus, if
the equation (\ref{eq44}) is satisfied by deformation functions
$H^m(x,\tilde t)$, it is also satisfied by functions
$H'^m(x,\tilde t)$, which are defined by formula (\ref{eq46}) with
arbitrary $\Lambda(x)^m{}_n$ from the group  $SO(m,n-m)^g$. All
such groups $\Gamma^{H'}_T$ specify the same (pseudo)Riemann
structure on $O$.

From geometric point of view, the field of (pseudo)orthonormal
vielbeins: $X'_m=\Lambda^{-1}(x)^n{}_mX_n$ changes during
transformations (\ref{eq46}). By the field of (pseudo)orthonormal
vielbeins $X_m$  from equation (\ref{eq44}) deformation functions
are \emph{uniquely} defined (let's recall that we assume
coordinates in $O$ to be fixed).

Let's point out that in our approach in (pseudo)Riemannian space
with respect to the field of (pseudo)orthonormal vielbeins  $X_m$
the rule of parallel transport of vectors to finite distance
$\tilde t(x)=x'-x=K(x,t(x))$ is uniquely defined. On the other
hand, in general case of curved space the result of parallel
transport depends upon the curve along which it is performed. So
there is a question to be asked: along which curve connecting
points  $x'$ and $x$ in the general case of curved
(pseudo)Riemannian space during performance of integral sequence
of infinitesimal transports (\ref{eq40}) do we get the result
which is given by formula (\ref{eq39})? We plan to study this
problem in our next publication \cite{8}.

In the first order with respect to $\tilde t$ equation
(\ref{eq44}) produces:
\begin{equation}\label{eq47}
\Gamma^\centerdot_{\mu\nu\sigma}+\Gamma^\centerdot_{\nu\mu\sigma}=\partial_\sigma
g_{\mu\nu},
\end{equation}
\noindent which, in view of symmetry of coefficients
$\Gamma^\sigma_{\mu\nu}$ in lower indices, gives formula:
\begin{equation}\label{eq48}
\Gamma^\centerdot_{\sigma\mu\nu}=\frac{1}{2}(\partial_\mu
g_{\nu\sigma}+\partial_\nu g_{\mu\sigma}-\partial_\sigma
g_{\mu\nu}),
\end{equation}
\noindent which, naturally, could be derived as the result of
formula (\ref{eq43}) in consideration of relation (\ref{eq36}).
Formula (\ref{eq48}) indicates that
$\Gamma^\centerdot_{\sigma\mu\nu}$  and $\Gamma^\sigma_{\mu\nu}$
in our case become \emph{Christoffel symbols of the $I^{st}$ and
$II^{nd}$ type} respectively.

In the second order with respect to $\tilde t$ it follows from
equation (\ref{eq44}) that:
\[\Delta^\centerdot_{\mu\nu\sigma\rho}+\Delta^\centerdot_{\nu\mu\sigma\rho}=
\partial^2_{\sigma\rho}g_{\mu\nu}-\Gamma^\centerdot_{\tau\mu\sigma}\Gamma^\tau_{\nu\rho}-
\Gamma^\centerdot_{\tau\nu\sigma}\Gamma^\tau_{\mu\rho},
\]
\noindent which, in view of symmetry of coefficients
$\Delta^\sigma_{\mu\nu\rho}$  in lower indices and relation
(\ref{eq47}), gives:
\[\Delta^\sigma_{\mu\nu\rho}=\frac{1}{3}(\partial_\rho
\Gamma^\sigma_{\mu\nu}+\partial_\nu
\Gamma^\sigma_{\mu\rho}+\partial_\mu
\Gamma^\sigma_{\nu\rho}+\Gamma^\sigma_{\tau\rho}\Gamma^\tau_{\mu\nu}+
\Gamma^\sigma_{\tau\nu}\Gamma^\tau_{\mu\rho}+\Gamma^\sigma_{\tau\mu}\Gamma^\tau_{\nu\rho}).
\]
\noindent Inserting this expression into formula (\ref{eq37}) and
considering formula (\ref{eq38}) we derive the expression
\[\rho^\sigma{}_{\mu\nu\rho}=\frac{1}{3}(R^\sigma{}_{\mu\nu\rho}+R^\sigma{}_{\nu\mu\rho}),
\]
\noindent the insertion of which into definition (\ref{eq27})
produces a well-known identity for curvature tensor:
\[R^\sigma{}_{\mu\nu\rho}+R^\sigma{}_{\nu\rho\mu}+R^\sigma{}_{\rho\mu\nu}=0.
\]

Thus we have proved

\textbf{Theorem 2.} \emph{Deformed group $\Gamma^H_T$  of
diffeomorphisms of chart $O$, produced with the help of
deformation, functions of which satisfy the equation (\ref{eq44}),
specifies by its action on $O$ a field of (pseudo)orthonormal
vielbeins and structure of (pseudo)Riemannian space. In
particular, coefficients of affine connection in coordinate basis
become equal to Christoffel symbols. The same structure of
(pseudo)Riemannian space is specified on $O$ by all the groups
$\Gamma^{H'}_T$, the deformation functions of which are connected
by transformations (\ref{eq46}) from gauge group  $SO(m,n-m)^g$.}

\emph{Arbitrary (pseudo)Riemannian structure on O can be specified
in such fashion.}

Through this theorem Klein's Erlangen Program is realized for
geometric structure of (pseudo) Riemannian space.

This work makes a group-theoretic description of geometric
structures of torsion-free affine connection and
(pseudo)Riemannian space locally within the bounds of a single
coordinate chart. This limitation can be lifted provided Lie
pseudogroups are studied.

The fact that relations derived from the conditions of existence
of certain groups have profound geometric meaning is another
confirmation of fundamentality of the ideas of Klein's Erlangen
Program in respect that geometry is defined completely by a group
of congruencies.


\end{document}